\documentclass[12pt]{amsart}

\usepackage[margin=3cm]{geometry}

\usepackage{amsmath,amssymb,amsbsy,amsfonts,amsthm,latexsym,
                        amsopn,amstext,amsxtra,euscript,amscd,mathrsfs,color,bm, cite}
                        \usepackage{ulem}

\begin{document}


\def\cA{{\mathcal A}}
\def\cB{{\mathcal B}}
\def\cC{{\mathcal C}}
\def\cD{{\mathcal D}}
\def\cE{{\mathcal E}}
\def\cF{{\mathcal F}}
\def\cG{{\mathcal G}}
\def\cH{{\mathcal H}}
\def\cI{{\mathcal I}}
\def\cJ{{\mathcal J}}
\def\cK{{\mathcal K}}
\def\cL{{\mathcal L}}
\def\cM{{\mathcal M}}
\def\cN{{\mathcal N}}
\def\cO{{\mathcal O}}
\def\cP{{\mathcal P}}
\def\cQ{{\mathcal Q}}
\def\cR{{\mathcal R}}
\def\cS{{\mathcal S}}
\def\cT{{\mathcal T}}
\def\cU{{\mathcal U}}
\def\cV{{\mathcal V}}
\def\cW{{\mathcal W}}
\def\cX{{\mathcal X}}
\def\cY{{\mathcal Y}}
\def\cZ{{\mathcal Z}}

\def\A{\mathbb{A}}
\def\B{\mathbf{B}}
\def \C{\mathbb{C}}
\def \F{\mathbb{F}}
\def \K{\mathbb{K}}

\def \Z{\mathbb{Z}}
\def \P{\mathbb{P}}
\def \R{\mathbb{R}}
\def \Q{\mathbb{Q}}
\def \N{\mathbb{N}}
\def \Z{\mathbb{Z}}

\newtheorem{theorem}{Theorem}
\newtheorem{lemma}[theorem]{Lemma}
\newtheorem{claim}[theorem]{Claim}
\newtheorem{cor}[theorem]{Corollary}

\title[Polynomial expanders  with many variables]{On  polynomial expanders with many variables}

\author[M.~Z.~Garaev]{M.~Z.~Garaev}
\address{Centro  de Ciencias Matem{\'a}ticas,  Universidad Nacional Aut\'onoma de
M{\'e}\-xico, C.P. 58089, Morelia, Michoac{\'a}n, M{\'e}xico}
\email{garaev@matmor.unam.mx}

\author[S. V. Konyagin] {S. V. Konyagin}
\address{Steklov Mathematical Institute of Russian Academy of Sciences, 8 Gubkina St., 119991 Moscow, Russia}
\email{konyagin@mi-ras.ru}

\maketitle

 \centerline{\it{In memory of Vera T. Sos}}

\bigskip

\begin{abstract} For a fixed integer $n\ge 2,$ we consider the homogeneous polynomial
$$
P(x_1, x_2, \ldots, x_{n+2})=\sum_{i=1}^{n} (x_2-x_1)^{i-1} x_1^{n-i} x_{i+2}.
$$
We prove that, for any finite set $A$ of complex numbers,
$$
\Bigl|\bigl\{P(x_1,x_2,\ldots,x_{n+2}): \, x_i\in A\bigr\}\Bigr|\gg |A|^{n}.
$$
The implicit constant in $\gg$ may depend only on $n.$
\end{abstract}

\bigskip


\section{\bf Introduction}

In what follows, $\R$ and $\C$ are the set of real and complex numbers. Given a finite set $X$ we denote by $|X|$ its cardinality. As usual, for sets $A$ and $B$, the
sum-set $A+B$ and the product-set $AB$ are defined by

$$
A+B=\{a+b:\, a\in A,\, b\in B)\},\quad AB=\{ab:\, a\in A,\, b\in B)\}.
$$

The sum-product estimate claims that there exists $c>0$ such that for any set $A\subset \R$ one has
$$
\max\{|A+A|, |AA|\} \gg |A|^{1+c}.
$$
The value of $c$ has been quantified in many papers. In the seminal work by Solymosi~\cite{Sol}, it is shown that any $c<1/3$ is admissible. Starting from the work of Konyagin and Shkredov~\cite{KonSh1}, the constant $1/3$ has been improved in a series of papers, with the most recent result, due to Cushman~\cite{Cush}, claiming that one can take any $c<\frac{1}{3} +\frac{10}{4407}.$
The famous conjecture of Erd\H{o}s and Szemer\'edi suggested that one could possibly take any fixed constant $c<1.$
Recently Bloom, Sawin, Schildkraut and Zhelezov~\cite{BSSZh} have surprisingly disproved this conjecture.
The authors indicate that they were inspired to revisit the possibility
of disproving the sum--product conjecture using number fields of large degree by the recent
OpenAI counterexample of the unit distance conjecture.
Currently it is not known what would be the candidate for the optimal value of $c.$

There are many versions of the sum-product problem, see~\cite{BaRo, BaRoZh, HanNetSt, TJ, MRSh, RoRu, RoWa, RoRuSh, RuSo}. Here, we record several of them with quadratic and superquadratic growth.

The results of Jones~\cite{TJ}: for any $A\subset\R,\, |A| \ge 2,$
$$
\Bigl|\Bigl\{\frac{a(b-c)}{c(a-b)}:\, a,b,c\in A\Bigr\}\Bigr|\gg \frac{|A|^2}{\log |A|}.
$$
and
$$
\Bigl|\Bigl\{\frac{(a-b)(c-d)}{(b-c)(a-d)}: \, a,b,c,d\in A\Bigr\}\Bigr|\gg |A|^2.
$$
Balog and Roche-Newton~\cite{BaRo}: for any $A\subset\R,\,|A| \ge 2,$
$$
\Bigl|\Bigl\{\frac{ab+c}{ad+ec}: a,b,c,d, e\in A\Bigr\}\Bigl|\gg\frac{|A|^{2+\frac{1}{8}}}{\log |A|}.
$$
This estimate, in particular, gives a superquadratic five variable expander with growth rate of $|A|^{2+\frac{1}{8}-o(1)}.$

Roche-Newton and Warren~\cite{RoWa}: for any $A\subset \R,$
$$
\Bigl|\Bigl\{\frac{ab-cd}{a-c}: a,b,c,d\in A\Bigr\}\Bigl|
\gg |A|^{2+ \frac{1}{14}}.
$$

Hanson, Roche-Newton and Steven~\cite{HanNetSt}: for any $A\subset\R,$
$$
|ab+(c-a)d: a,b,c,d \in A|\gg |A|^2.
$$

A particular case of the breakthrough paper by Guth and Katz~\cite{GuKa} on distinct distance problem implies that for any $A\subset \R,$
$$
|\{(a-b)^2+(c-d)^2: \, a,b,c,d\in A\}|\ge |A|^{2-o(1)}.
$$

\bigskip

In the present paper, we are interested in polynomial expanders with many variables over sets of complex numbers.

In what follows, the implied constants in $\ll$ and
$\gg$ symbols may depend only on $n.$

\section{Statement of our result}

Our result is as follows.

\begin{theorem}
\label{thm:Main}
Let $n\ge 2$ be a fixed integer and let
$$
P(x_1, x_2, \ldots, x_{n+2})=\sum_{i=1}^{n} (x_2-x_1)^{i-1} x_1^{n-i} x_{i+2}.
$$
Then, for any finite set $A\subset \C,$ we have
$$
\Bigl|\bigl\{P(x_1,x_2,\ldots,x_{n+2}): \, x_i\in A\bigr\}\Bigr|\gg |A|^{n}.
$$
\end{theorem}

In other words,  our homogeneous polynomial of degree $n$ and in $n+2$ variables expands subsets of $\C$ with the $n$-th power growth rate. The example
$A={[1,N]\cap\Z}$ shows that the lower bound of our  estimate is sharp. It would be very interesting to find a polynomial with $k<n+2$ variables which would expand sets $A\subset \C$ with the same rate of growth $|A|^n.$ Clearly, such a number $k$ can not be smaller than $n.$

For $n=2$ and $A\subset\C$ we get
$$
\Bigl|\{x_1x_3+(x_2-x_1)x_4:  \, x_1,x_2,x_3,x_4\in A\}\Bigr|\gg |A|^2.
$$
This extends the aforementioned result of~\cite{HanNetSt} from the real settings to complex numbers.

For $n=3$ we have the following $5$-variable polynomial, which expands sets  $A\subset \C$ with the cubic growth rate:
$$
\Bigl|\bigl\{x_1^2x_3+(x_2-x_1)x_1x_4+(x_2-x_1)^2x_5: x_i\in A\bigr\}\Bigr|\gg |A|^3.
$$

\section{Proof of Theorem~\ref{thm:Main}}

We can assume that $|a|\ge 2$ for any $a\in A$ and that $|A|$ is sufficiently large in terms of $n.$
Furthermore, we can assume that
\begin{equation}
\label{eqn:arg x/y}
|\arg (x/y)| <\pi/3,\quad \text {for any} \,\, x,y\in A.
\end{equation}
Indeed, we can split the plain into $6$ angle sectors with the angle $\pi/3$ and the common vertex at the origin.  One of these sectors contain at least $|A|/6$ elements of
$A.$ We can confine ourselves to this subset of $A,$ which without loss of generality we again denote  by $A.$ Thus, we can assume that~\eqref{eqn:arg x/y} holds.

Take a large number $M$ depending on $n.$

Let $a,b\in A,\, a\not=b$ be chosen such that
$$
\rho=\Bigl|1-\frac{b}{a}\Bigr| = \min_{x,y\in A,\, x\not=y}\,\Bigl|1-\frac{x}{y}\Bigr|.
$$
Note that
\begin{equation}
\label{eqn: b le a}
|b| \le |a|.
\end{equation}
Indeed, by the definition of $\rho$ we have
$$
\Bigl|1-\frac{b}{a}\Bigr| \le \Bigl|1-\frac{a}{b}\Bigr|.
$$
On the other hand,
$$
\Bigl|1-\frac{b}{a}\Bigr|  = \frac{|b|}{|a|}\Bigl|1-\frac{a}{b}\Bigr|,
$$
and (\ref{eqn: b le a}) follows.

Next, note that
$$
0<\rho<1.
$$
Indeed, let $z = b/a$. By (\ref{eqn: b le a}), $|z|\le 1.$ For $\theta =\arg (b/a)$ we have $|\theta|<\pi/3.$ Then
$$
\rho\le |1-z|=\sqrt{1-2|z|\cos\theta+|z|^2}< \sqrt{(1-|z|+|z|^2}\le 1.
$$

To prove our theorem, we consider two cases.

\bigskip

{\bf Case 1.}  Assume that $\rho \ge 1/(4 M)$.

\bigskip

Consider the ring sector $r \le |z| \le 2r.$ Let $G$ be the elements of $A$ contained in this sector. 
Then for any distinct elements $z_1, z_2\in G,$
$$
|z_1-z_2|\ge \rho |z_2|\ge \rho r.
$$
Consider $|G|$ open discs with centers at the points of $G$ and radius $0.5\rho r$.
 All these  discs are pairwise disjoint and are contained inside the disc $|z|\le 2r+0.5\rho r.$ Hence,
$$
|G|\pi(0.5\rho r)^2<\pi (2r+0.5\rho r)^2,
$$
whence
$$
|G|<\frac{(\rho+4)^2}{\rho^2}<400 M^2.
$$
Thus any ring sector $r \le |z| \le 2r$ contains less than $400 M^2$ elements of
$A$.

Split the region $|z|\ge 2$ into the ring sectors of the form
$$
2^k\le |z|<2^{k+1}.
$$
Then $A$ is contained in the union of $4$ regions of the form
$$
R_\ell= \bigcup_{k=\ell\bmod 4}\{2^{k}\le |z|< 2^{k+1}\},\quad \ell =0,1,2,3.
$$
One of these regions contains at least $|A|/4$ elements of $A.$ Denote this region by  $R_{\ell}.$ Since each of the ring sector $2^{k}\le |z|< 2^{k+1}$ 
contains at most $400 M^2$ elements of $A,$ we can construct $A_1\subset A\cap R_\ell$ such that
$$
|A_1|\gg |A|
$$
and  $A_1 = \{ x_1, x_2,\dots\}$ with $|x_j| \ge 8|x_{j+1}|$.

The number of $n$-tuples $(a_1, a_2,\ldots, a_n)$ with
\begin{equation}
\label{eqn:small rho a1>a2>an}
a_1,\ldots, a_n \in A_1,\quad |a_1|> |a_2| > \ldots > |a_n|,
\end{equation}
is clearly $\gg |A_1|^n\gg |A|^n.$ Thus, it suffices to prove that all the numbers
$$
a^{n-1}a_1 + a^{n-2}(b-a)a_2+\ldots + (b-a)^{n-1}a_n,
$$
with $a_1, a_2,\ldots, a_n$ satisfying~\eqref{eqn:small rho a1>a2>an}, are pairwise distinct. To prove this, assume that
$b_1, b_2,\ldots, b_n\in A_1$ are such that
$$
|b_1|> |b_2| > \ldots > |b_n|
$$
and
$$
\sum_{i=1}^{n}a^{n-i}(b-a)^{i-1}a_i = \sum_{i=1}^{n}a^{n-i}(b-a)^{i-1}b_i.
$$
It suffices to prove that $a_i=b_i$ for all $1\le i\le n.$ Assume contrary, let $j$ be the smallest index such that $a_j\not =b_j.$
Then cancelling out by $(b-a)^{j-1},$ we get
$$
\sum_{i=j}^n a^{n-i}(b-a)^{i-j}a_i = \sum_{i=j}^n a^{n-i}(b-a)^{i-j}b_i.
$$
By the construction of the set $A_1,$ we can assume that $|a_j|>|b_j|,$ which imples that $|a_j|\ge 8|b_j|$. Separating the terms corresponding to $i=j,$ we get
\begin{equation}
\label{eqn:separate term i=j}
a^{n-j}(a_j-b_j)= \sum_{i=j+1}^n a^{n-i}(b-a)^{i-j}(b_i-a_i).
\end{equation}
Note that
\begin{equation}
\label{eqn:estimate LHS}
|a^{n-j}(a_j-b_j)|\ge |a|^{n-j}(|a_j|-|b_j|)\ge \frac{7|a|^{n-j}|a_j|}{8}.
\end{equation}
Furthermore, for $i\ge j+1$ we have that
$$
|b_i-a_i|\le |b_i|+|a_i|\le \frac{|b_j|}{8^{i-j}} +\frac{|a_j|}{8^{i-j}}\le \frac{2|a_j|}{8^{i-j}}.
$$
Hence,
$$
|a^{n-i}(b-a)^{i-j}(b_i-a_i)|\le |a|^{n-i}\rho^{i-j}|a|^{i-j}(|b_i|+|a_i|)\le \frac{2|a|^{n-j}|a_j|}{8^{i-j}}.
$$
Therefore,
$$
\Bigl|\sum_{i=j+1}^n a^{n-i}(b-a)^{i-j}(b_i-a_i)\Bigr|\le 2|a|^{n-j}|a_j|\sum_{i=j+1}^{n}\frac{1}{8^{i-j}}<\frac{2|a|^{n-j}|a_j|}{7}.
$$
This and~\eqref{eqn:estimate LHS} implies that
$$
\Bigl|\sum_{i=j+1}^n a^{n-i}(b-a)^{i-j}(b_i-a_i)\Bigr|< |a^{n-j}(a_j-b_j)|,
$$
which contradicts to~\eqref{eqn:separate term i=j}.

\bigskip

{\bf Case 2.} Let now  $\rho < 1/(4 M)$.

\begin{lemma}
\label{lem:small rho} If $\rho < 1/(4 M)$, then there is a set $A_1 \subset A$ such that $|A_1| \gg |A|$
and for any distinct points $a',b'\in A_1$ we have
$$
\Bigl|1 - \frac{b'}{a'}\Bigr| \ge M\rho.
$$
\end{lemma}

\begin{proof}

Take any element $z_1\in A.$ Let $G_1$ be the set of elements $y\in A$ with
$$
|z_1-y|\le 2M \rho y.
$$
Since $2M \rho < 0.5$, we have
$$
0.5|z_1|<|y|<2 |z_1|.
$$
In particular,
$$
|z_1-y|\le 4M\rho |z_1|.
$$
It follows that $G_1$ is contained in the disc with the center at $z_1$ and radius $4M\rho|z_1|.$

Consider the collection of $|G_1|$ discs with centers at points of $G_1$ and radius $0.1\rho |z_1|.$
These discs are pairwise disjoint. Indeed, if two of these discs, with centers at $y_1$ and $y_2$ have non-empty intersection, then
$$
|y_1-y_2|\le 0.2\rho|z_1|\le  0.4\rho |y_2|,
$$
implying
$$
\Bigl|1-\frac{y_1}{y_2}\Bigr|\le 0.4\rho.
$$
This contradicts the definition of $\rho.$ Hence, all these discs are pairwise disjoint. Furthermore, all them are contained inside the disc with center at $z_1$ and radius $4M\rho|z_1|+0.1\rho|z_1|.$ Hence,
$$
|G_1|\pi (0.1\rho |z_1|)^2<\pi (4M\rho|z_1|+0.1\rho z_1)^2.
$$
It follows that
$$
|G_1|<2000 M^2.
$$


We now construct the set $A_1$ as follows. We include the element $z_1$ into the set $A_1$ and
from the set $A$ remove the element $z_1$ and all the elements $y\in A$ satisfying $|z_1-y|\le 2M\rho y.$ Since $|G_1|<2000M^2$,  this way we removed at most $2000 M^2$ elements from $A$ but added $1$ element to $A_1.$

We repeat the process with the new  remaining set $A$ and continue the process until we construct the set $A_1$ and all elements of $A$ are exhausted. Clearly, we end up with a set $A_1$ satisfying
$$
|A_1|\ge \frac{|A|}{2000 M^2}\gg |A|,
$$
and such that for any two elements $z_1, z_2\in A_1$ one of the two inequalities hold:
$$
|z_1-z_2| > 2M\rho |z_1|,\quad {\text or} \quad  |z_1-z_2| > 2M\rho |z_2|.
$$
We shall prove that in either case we simultaneously have
$$
|z_1-z_2| > M\rho |z_1|,\quad \text {and} \quad  |z_1-z_2| > M\rho |z_2|.
$$
Let, say, $|z_1-z_2| > 2M\rho |z_1|.$ We need to prove that $|z_1-z_2| > M\rho  |z_2|.$ We have
$$
z_2 =z_1-\Delta M\rho z_1,\quad |\Delta|>2.
$$
Then,
$$
\Bigl|\frac{z_1}{z_2}-1\Bigr|=\Bigl|\frac{1}{1-\Delta M \rho }-1\Bigr| =\frac{|\Delta M\rho|}{|1-\Delta M\rho|}\ge \frac{|\Delta|}{1+|\Delta| M \rho}M\rho > M\rho.
$$
Here we used the fact that $\rho< 1(4 M)$ and hence
$$
1+|\Delta|M\rho\le \frac{|\Delta|}{2}+\frac{|\Delta|}{4}<|\Delta|
$$
\end{proof}

In order to deal with the Case 2, we choose $A_1\subset A$ as in lemma. It suffices to prove that all numbers of the form
$$
a^{n-1}a_1 + a^{n-2}(b-a)a_2+\ldots + (b-a)^{n-1}a_n,
$$
with pairwise distinct $a_1, a_2,\ldots, a_n\in A_1$ such that $|a_1|\ge |a_2|\ge \ldots |a_n|,$ are pairwise distinct. To prove this, assume that
$b_1, b_2,\ldots, b_n\in A_1$ are such that
$$
|b_1|\ge |b_2| \ge \ldots \ge |b_n|
$$
and
$$
\sum_{i=1}^{n}a^{n-i}(b-a)^{i-1}a_i = \sum_{i=1}^{n}a^{n-i}(b-a)^{i-1}b_i.
$$
It suffices to prove that $a_i=b_i$ for all $1\le i\le n.$ Assume contrary, let $j$ be the smallest index such $a_j\not =b_j.$ We can assume that $|a_j|\ge |b_j|.$ Then,
\begin{equation}
\label{eqn:final equality}
a^{n-j}(a_j-b_j)= \sum_{i=j+1}^n a^{n-i}(b-a)^{i-j}(b_i-a_i).
\end{equation}
We estimate the absolute value of each term on the right-hand-side as follows. Since $i\ge j+1$ and $\rho <1,$ we have
\begin{equation*}
\begin{split}
|a^{n-i}(b-a)^{i-j}(b_i-a_i)|&\le |a|^{n-i}(\rho|a|)^{i-j}(|b_i|+|a_i|)\\ &\le \rho|a|^{n-j}(|b_j|+|a_j|)\\ &\le 2\rho|a|^{n-j}|a_j|.
\end{split}
\end{equation*}
Thus,
$$
\Bigl|\sum_{i=j+1}^n a^{n-i}(b-a)^{i-j}(b_i-a_i)\Bigr|\le 2n\rho |a|^{n-j}|a_j|.
$$
On the other hand, by Lemma~\ref{lem:small rho}, we have
$$
|a^{n-j}(a_j-b_j)|\ge |a|^{n-j}M\rho |a_j|.
$$

Hence, for sufficiently large $M$ we get
$$
|a^{n-j}(a_j-b_j)|> \Bigl|\sum_{i=j+1}^n a^{n-i}(b-a)^{i-j}(b_i-a_i)\Bigr|,
$$
which contradicts \eqref{eqn:final equality}. This contradiction shows that $a_i=b_i$ for all $1\le i\le n$ and finishes the proof of the Case 2 of our theorem.

\bigskip

{\bf Acknowledgement.}

The work of the second author was performed at the Steklov
International Mathematical Center and supported by the Ministry
of Science and Higher Education of the Russian Federation
(agreement no. 075-15-2022-265).

\end{document}